\documentclass[english]{amsart}

\usepackage{babel}
\usepackage{amstext}
\usepackage{amsmath}
\usepackage{amsfonts}
\usepackage{latexsym}
\usepackage{ifthen}
\usepackage{xypic}
\xyoption{all}
\newtheorem{lemma1}[equation]{}

\newenvironment{lemma}{\begin{lemma1}{\bf Lemma.}}{\end{lemma1}}

\newenvironment{example}{\begin{lemma1}{\bf Example.}\rm}{\end{lemma1}}

\newenvironment{theorem}{\begin{lemma1}{\bf Theorem.}}{\end{lemma1}}

\newenvironment{proposition}{\begin{lemma1}{\bf Proposition.}}{\end{lemma1}}
\newenvironment{corollary}{\begin{lemma1}{\bf Corollary.}}{\end{lemma1}}

\newenvironment{remark}{\begin{lemma1}{\bf Remark.}\rm}{\end{lemma1}}

\newenvironment{setup}{\begin{lemma1}{\bf Setup.}}{\end{lemma1}}

\newenvironment{conjecture}{\begin {lemma1}{\bf Conjecture.}}{\end{lemma1}}

\makeatletter
\ifnum\@ptsize=0  \fi \ifnum\@ptsize=2  \fi 

\newcommand\sE{{\mathcal E}}
\newcommand\sF{{\mathcal F}}

\newcommand\sI{{\mathcal I}}

\newcommand\sL{{\mathcal L}}
\newcommand\sO{{\mathcal O}}
\newcommand\sS{{\mathcal S}}

\newcommand\bQ{{\mathbb Q}}
\newcommand\bN{{\mathbb N}}

\title {Kodaira dimension of subvarieties II} 
\author{Thomas Peternell}

\begin{document}

\maketitle

\tableofcontents

\section{Introduction}
This paper continues the study of non-general type subvarieties begun in [PSS99]. Consider for simplicity a submanifold $A$ of a projective
manifold $X$ with ample normal bundle $N_A$ (but we will also consider also subvarieties with weaker positivity assumptions on the normal bundle). 
In [PSS99] it was shown that if $A$ is not of general type, then the Kodaira dimension
$\kappa (X) = - \infty.$ Conjecturally any manifold of negative Kodaira dimension is uniruled, i.e. covered by rational curves. 
Here we establish that in the above situation. The new ingredient is the paper [BDPP04] where it is shown that a manifold is uniruled
if its canonical bundle $K_X$ is not pseudo-effective, i.e. its Chern class is not contained in the closure of the cone generated
by effective classes. Alternatively, $K_X$ is not pseudo-effective, if $K_X$ does not admit a possibly singular metric with non-negative
curvature (current). \\
Once $X$ is uniruled, one can form the rational quotient $f: X \dasharrow W,$ see [Ca81],[Ko96]. 
This is just a rational map, but has some good properties, see
(2.1). Now by a fundamental result of Graber-Harris-Starr [GHS03], $W$ is no longer uniruled, hence again by [BDPP04]
$K_W$ is pseudo-effective ($W$ can be taken smooth).
Coming back to our submanifold $A,$ this can be used to show that $f \vert A$ is
onto $W$; thus we obtain a bound for $\dim W$ and often also for $\kappa (W).$ It allows also to relate $\pi_1(A)$ and $\pi_1(X).$  Finally
we give a cohomological criterion for a projective manifold to be rationally connected.\\
For interesting discussions I would like to thank Fr\'ed\'eric Campana; I also thank the referee for important remarks. The paper was
written within the framework of the priority program ``Global Methods in Complex Geometry'' of the Deutsche Forschungsgemeinschaft. 

\section{Preliminaries}

{\bf (2.1).} Let $X$ be a uniruled projective manifold. By identifying sufficiently general points which can be joined by chains
of rational curves, one obtains the so-called rational quotient (e.g. [Ca81], [CP99], [Ko96]) 
$$ f: X \dasharrow W $$
which has the remarkable property to be almost holomorphic, ie. the locus of indeterminacy does not project onto $W.$ 
By [GHS03], $W$ does not carry any covering family of rational curves, ie. $W$ is not uniruled, and consequently the canonical
bundle $K_W$ is not pseudo-effective (if $W$ is smooth, [BDPP04]). \\
If $\dim W = 0,$ then $X$ is {\it rationally connected}.

\vskip .2cm \noindent 
{\bf (2.2).} In [Ca95], Campana introduced a refined Kodaira dimension, which potentially differs from the usual notion only in the
case that $X$ is uniruled. Namely, he considers subsheaves $\sF  \subset \Omega^p_X $ of an arbitrary rank $r > 0$ for any $p$ and defines $$ \kappa^+(X) $$
to be the maximum over all $\kappa (\det \sF)$, where $\det \sF = (\bigwedge^r \sF)^{**}.$ 
Conjecturally $\kappa^+(X) = \kappa (X)$ if $X$ is not uniruled; and $\kappa (X) = \kappa (W),$ when $X$ is uniruled with rational quotient $W$ -
here we let $\kappa (W) = -\infty$ if $\dim W = 0.$  See [Ca95],[CP99] for details. In [CP04] it is proved that $\kappa^+(X) = \dim X$ implies 
$\kappa (X) = \dim X.$

\vskip .2cm \noindent 
{\bf (2.3)} (1) We say that a vector bundle $E$ is $\bQ-${\it effective} [PSS99], if some symmetric power $S^k(E)$ is generically
spanned by global sections. \\
(2) $E$ is {\it almost ample}, if $E$ is nef and $E \vert C$ is ample on every curve not contained in some countable union of proper
subvarieties. \\
(3) $E$ is {\it generically ample}, if $E$ is $\bQ-$effective and if $E \vert C$ is ample for all curves not contained in a 
countable union of proper subvarieties. \\
(4) $E$ is {\it big} if for all torsion free quotients $E \to \sS \to 0$ of rank 1 the determinant $\det \sS := \sS^{**} $ is big. 

\vskip .2cm \noindent 
As an example, consider a birational map $\pi: \hat X \dasharrow X$ from a projective manifold $\hat X.$ Consider a vector bundle $E$
of the form $E = \pi_*(\hat E)^{**}$ with an ample vector bundle $\hat E$ on $\hat X.$ Then $E$ is $\bQ-$effective, big and generically
ample. \\
Notice also that an almost ample vector bundle might not be big. 

\section{A criterion for uniruledness}

For the entire section we fix a projective manifold and a compact submanifold $A \subset X$ whose normal bundle $N_A = N_{A/X}$ is supposed to be
$\bQ-$effective. One could allow certain singularities
(see remark 3.7(4)) but for sakes of simplicity we only consider the smooth case.

\begin{lemma} Let $L$ be a pseudo-effective line bundle on $X$. Then the restriction $L_A$ is again pseudo-effective. 
\end{lemma}

\begin{proof} Fix an ample line bundle $H$ on $X$. Since $L$ is pseudo-effective, we find for any $m \gg 0$ a positive number $t$ such that
$$ H^0(X,t(mL+H)) \ne 0.$$ 
Then from [PSS99,2.1] we deduce the existence of some number $k$ such that 
$$ H^0(S^kN^*_A \otimes t(mL_A+H_A)) \ne 0.$$
By assumption there exists $l$ such that $S^lN_A$ is generically spanned, hence $S^lN_A^* \subset \sO_A^q.$
Putting things together, we obtain a number $t'$ such that
$$H^0(t'(mL_A+H_A)) \ne 0.$$ 
Thus $L_A$ is pseudo-effective.
\end{proof}

As application we prove

\begin{theorem} If $N_A$ is $\bQ-$effective and $\det N_A$ is big and $A$ not of general type, then $X$ is uniruled.
\end{theorem}

\begin{proof} We are going to show that $K_X$ is not pseudo-effective; then the uniruledness follows from [BDPP04]. Assume to the contrary that
$K_X$ is pseudo-effective. Then by the previous lemma $K_X \vert A$ is pseudo-effective. Thus by adjunction, $K_A$ is the sum of a 
pseudo-effective and a big divisor, hence big, contradicting our assumption.
\end{proof}

\begin{remark} In [PSS99] it is shown that if $N_A$ is $\bQ-$effective and $\kappa (\det N_A) > \kappa (A),$ then $\kappa (X) = - \infty.$
So we would expect that in Theorem 3.2 that it is sufficient to assume $\kappa (\det N_A) > \kappa (A).$ 
>From our discussion in the proof we see that things come down to the following. Suppose that 
$$ K_A = P + E$$
with $P$ pseudo-effective and $E$ effective, then we should have $\kappa (A) \geq \kappa(E).$ This is discussed in [CP04] and proved in special cases.
It is not unreasonable to expect that the assumption on the $\bQ-$effectivity can be substituted by assuming that $N_A$ is nef. 
\end{remark}

\begin{corollary} Let $X$ be a projective manifold, $A \subset X$ a submanifold with ample normal bundle. If $A$ is not of general type, then
$X$ is uniruled. 
\end{corollary} 

\begin{setup}  {\rm  Returning to submanifolds $A \subset X$ with $\bQ-$effective normal bundle $N_A$ and big determinant $\det N_A$  
and assuming again $A$ not of
general type, we now know that $X$ is uniruled. So let 
$$ f: X \dasharrow W$$
be ``the'' rational quotient; we choose $W$ to be smooth. Since $W$ is not uniruled, $K_W$ is pseudo-effective and if $p = \dim W$, there exists a
locally free pseudo-effective sheaf $\sF \subset \Omega^p_X $ of rank 1. To be more precise, consider a birational map $\pi: \hat X \to X $ from a
projective manifold $\hat X$ such that the induced map $\hat f: \hat X \to W$ is holomorphic. Then we have an injective map
$ \hat f^*(\Omega_W^p) \to \Omega^p_{\hat X}$ and we set
$$ \sF = (\pi_*(\hat f^*(K_W)))^{**}.$$ 
Then $\sF \subset \Omega^p_X$ since this evidently holds outside a set of codimension at least $2$ and $\sF$ is clearly pseudo-effective. 
Moreover, if $A \subset X$ with $N_A$ $\bQ-$effective is a subvariety then Lemma 3.1 shows that $\sF_A$ is again
pseudo-effective.} 
\end{setup} 
If $A$ does not meet the set of indeterminacies of $f,$ then the {\it almost ampleness} of $N_A$ forces $f \vert A$ to be onto as soon as 
$\dim A \geq \dim W.$ This actually also holds if $A$ does meet the indeterminacy locus:

\begin{lemma} Let $X$ and $W$ be projective manifolds, $f: X \dasharrow W$ be almost holomorphic. Let $A \subset X$ be a submanifold or a local complete
intersection with almost
ample normal bundle and assume $\dim A \geq  \dim W.$ Then $A$ meets the general fiber of $f$, in particular $A$ is not contained in the
indeterminacy locus of $f$, and $f: A \dasharrow W$ is dominant. 
\end{lemma}

\begin{proof} If $N_A$ is ample, this is a special case
of [FL82],Theorem 1 or Corollary 1, appplied to a general fiber or a general hyperplane section of a general fiber. See also
[Fu84,12.2.4]. \\
The more general case that $N_A$ is almost ample can be done by an easy adaption of the proof of the ample case. 
\end{proof}

It is not clear whether the lemma is still true when $N_A$ is only $\bQ-$effective (or generically ample). \\
In the next theorem we make use of the notations of (3.5). 

\begin{theorem}  Assume that $A$ is not of general type with $N_A$ big and $\bQ-$effective.
Then:  
\begin{enumerate} 
\item $\sF_A$ is a (pseudo-effective) subsheaf of $\Omega^p_A.$
\item $\dim A \geq \dim W$, and if $N_A$ is almost ample or if $A$ is not in the indeterminacy locus of $f,$ then $f \vert A$ is onto $W$.
\item If $A$ is uniruled with rational quotient $A \dasharrow B,$ then even $\dim B \geq \dim W.$ 
If $N_A$ is ample or if $A$ is not in the indeterminacy locus of $f,$ then 
$f \vert A$ is onto with positive-dimensional fibers; more generally there is a surjective meromorphic map $B \dasharrow W.$ 
\end{enumerate}
\end{theorem}

\begin{proof} Recall that indeed $X$ is uniruled and consider the rational quotient $f: X \dasharrow W.$ First we show how to derive (2),(3) and (4) from (1). \\
(2) is obviously a consequence of (1); for the second part use Lemma 3.6  and observe that since $\sF_A$ is a subsheaf of $\Omega^p_A$ via the 
canonical generically defined map 
$f^*(\Omega^p_W) \to \Omega^p_A,$ the restricted map $f \vert A$ must be onto. \\
Concerning (3): if $A$ is uniruled, then we cannot have $\dim A = p,$ because then $K_A$ would contain a pseudo-effective line bundle,
i.e. would be pseudo-effective. More generally, consider a general smooth fiber$A_b$ of the almost holomorphic map $g: A \dasharrow B.$
Then $\Omega^p_A \vert A_b$ has a filtration by terms 
$$\Omega^i_{A_b} \otimes \bigwedge^j N^*_{A_b}. $$
Since $A_b$ is rationally connected, $\sF_A$ cannot be in a term $\Omega^i_{A_b} \otimes \bigwedge^j N^*_{A_b}$ with $i > 0,$ hence 
we must have $\sF_A \subset g^*(\Omega^p_B)$, at least generically. But then $\dim B \geq p.$ 
If $A$ is not contained in the indeterminacy locus of $f,$ then it is clear that the fibers $A_b$ must be contracted by $f$ (consider general
rational curves in $A_b;$ their deformations cover all of $X$). Hence $X \dasharrow W$ factors over $B \dasharrow W.$ 
\vskip .2cm \noindent 
It remains to prove (1).
By taking $\bigwedge^p$ of the exact sequence
$$ 0 \to N^*_A \to \Omega^1_X \vert A \to \Omega^1_A \to 0 $$
we obtain a filtration of $\Omega^p_X \vert A$ by the terms
$$ \Omega^j_A \otimes S^{i}N_A^* $$
with $i+j = p.$ The map $\sF \to \Omega^p_X $ corresponds to a non-zero section $s \in H^0(X,\Omega^p_X \otimes \sF^*).$ 
If $s \vert A = 0,$ choose $k$ maximal such that 
$$s \in H^0(X,\Omega^p_X \otimes \sF^* \otimes \sI_A^k).$$
Allowing $k = 0$ if $s \vert A \ne 0,$ we always get a non-zero map
$$ S^kN_A \otimes \sF_A \to \Omega^p_X \vert A.$$ 
This yields a non-zero map 
$$ S^kN_A \otimes \sF_A \to  \Omega^j_A \otimes S^{i}N_A^*$$
and therefore a non-zero map
$$ S^kN_A \otimes S^{i}N_A \otimes \sF_A \to \Omega^j_A.$$ 
Since $\sF_A$ is pseudo-effective (see 3.5) and since $N_A$ is assumed to be big, the determinant $\sL$ of the image of this map is big.  
If now $i > 0$ or $k > 0,$ then by [CP04], we obtain $\kappa (A) = \dim A,$ contrary to our assumption. 
Hence $\sF_A \subset \Omega^p_A.$

\end{proof}

\begin{corollary} Assume that $A$ is not of general type and $N_A$ to be $\bQ-$effective. Then the conclusions of Theorem 3.7 hold,
if one of the following conditions is satisfied.
\begin{enumerate}
\item $N_A$ is ample.
\item $\det N_A$ big; $A$ is not contained in the indeterminacy locus of the rational quotient; 
$W$ is of general type and $f \vert A$ is generically finite.
\item $A$ is uniruled, $N_A$ is generically ample.
\item $A$ is uniruled, $\kappa (W) > 0$.
\item For any map 
$$ S^kN_A \otimes \sF_A \to \Omega^j_A,$$
the determinant of the image is big.
\end{enumerate}
\end{corollary}  

\begin{proof}  (1) is clear. For (2) we only need to remark that $\sF_A$ is big;  then the proof of (3.7) still works. 
The bigness of $\sF_A$ follows from the fact that $f \vert A: A \dasharrow W$ is generically finite, so that it suffices to show
that $K_W \vert \overline {f(A)}$ is big. This in turn follows from the bigness of $K_W$, the fact that the normal sheaf of 
$\overline {f(A)} \subset W$ is $\bQ-$effective and [PS99,4.2]. \\
For (3) and (4),  consider a general rational curve $C$ in $A$ and first observe that the deformations of $C$ fill up $X$ so that
$X$ is uniruled. Then 
consider the restricted morphism
$$ \kappa_C: S^kN_A \vert C \otimes \sF_C \to \Omega^j_A \vert C \otimes S^{i}N_A \vert C.$$ 
Since $\sF_A$ is pseduo-effective (3.5) and since  $N_A \vert C$ is ample of $N_A$ is generically ample, $\kappa_C$ must vanish, the dual of $\Omega^j_A \vert C$ being nef. 
In case $N_A$ is just $\bQ-$effective, we can only conclude that either $\kappa_C = 0$ for general $C$ or does not have any zeroes. 
Hence if the full map $\kappa \ne 0$, then it can vanish only in codimension $2$. Then
consider a general curve $B$ so that $\kappa_B$ does not have zeroes. This is only possibly when $\sF_B$ is not ample which implies $\sF \equiv 0,$
hence $\kappa (W) = 0.$ \\
(5) finally is a direct consequence from the proof of (3.7)(1). 
\end{proof} 

The critical point in the proof of Theorem 3.7 where we use the bigness of $N_A$ is the investigation of a line bundle $\sL \subset \bigwedge^r \Omega^j_A.$ If say $N_A$ is 
generically ample, then this line bundle has the property that 
$$ L \cdot C > 0 $$
for all movable curves $C.$ Recall that a curve is movable, if the deformations of $C$ fill up the whole variety. 
In general, it is not true that a line bundle with this property is big. However the following should be true.

\begin{conjecture} Let $A$ be a projective manifold, $\sL \subset \bigwedge^r \Omega^j_A$ a line bundle such that $\sL \cdot C > 0$ for all 
movable curves. Then $\sL$ is big, in particular $A$ is of general type. 
\end{conjecture} 

Even if $L = K_A$, this conjecture is unknown in dimension at least $4$. However:

\begin{proposition} 
Assume that $A$ admits a good minimal model $A'$, i.e., some multiple $mK_{A'}$ is spanned. Then Conjecture 3.9 holds in the weak form that
$A$ is of general type.
\end{proposition} 

\begin{proof} It is easily checked that the line bundle $\sL$ induces a Weil divisor $\sL'$ on $A$ which is $\bQ-$Cartier such that
$$ \sL' \subset \bigwedge^r \Omega^j_{A'}.$$
This last sheaf is by definition the reflexive extension of the corresponding sheaf on the smooth part. 
Notice that the existence of $\sL$ forces $A$ to be non-uniruled. 
Consider the exact sequence
$$ 0  \to \sL' \to \bigwedge^r \Omega^j_{A'} \to Q \to 0.$$
Then $\det Q$ is generically nef by [CP04,1.5]. Taking determinants, we obtain 
$$ NK_{A'} = \sL' + \det Q$$  
for a suitable positive number $N.$
Thus $K_{A'} \cdot C > 0$ for all movable curves. Since some multiple of $K_{A'}$ is spanned, $A$ must be of general type. 
\end{proof} 

\begin{corollary} Assume $A \subset X$ not of general type and $\dim A  \leq 3.$ If $N_A$ is $\bQ-$effective and generically ample with
$\det N_A$ big, or, then Theorem 3.7 holds. 
\end{corollary}

\begin{remark} (1) The proof of (3.7) shows that $ \kappa (W) \leq \kappa^+(A)$ (with the assumption that $A $ is not in the indeterminacy
locus of $f$ in case $N_A$ is not almost ample). \\
(2) Via the conjectural equality $\kappa^+(X) = \kappa (W),$ one would arrive at 
$$ \kappa^+(X) \leq {\rm min}(\dim A,\kappa^+(A)).$$
(3) There should be a version of (3.7) in the case $\kappa (\det N_A) > \kappa (A).$ \\
(4) Everything in this section works also if $X$ and $A$ have canonical singularities and ${\rm codim}_A (A \cap {\rm Sing}(X)) \geq 2.$  

\end{remark}

\begin{example} {\rm Let $X$ be a projective manifold and $C \subset X$ an elliptic curve with ample normal bundle. Then $\dim W \leq 1,$
i.e. $\kappa^+(X) \leq 0.$ If $X$ is not rationally connected, i.e. $\dim W = 1,$ then the rational quotient $X \to W$ is holomorphic
and $C$ is multi-section. In particular $W$ is elliptic and $q(X) = 1.$ }
\end{example} 

\vskip .2cm \noindent
We will generalize that in the following section.

\setcounter{equation}{0}
\section{Computing invariants} 

In this section we compare the fundamental group and the spaces of $p-$forms for $A$ and $X$, when the normal bundle $N_A$ is almost ample and $\bQ-$effective. 
If $N_A$ is ample, a theorem of Napier-Ramachandran [NR98] says that, given a submanifold $A \subset X$ with ample
normal bundle, then the image of $\pi_1(A) \subset \pi_1(X)$ has finite index (instead of ampleness it actually suffices to assume 
finiteness of formal cohomology in degree $0$ along $A$ of locally free sheaves; but the relation of this property to almost ampleness is unclear).  
We find Napier-Ramachandran's theorem again if $A$ is not of general type, but we can also
weaken the ampleness assumption. Furthermore we can deal with holomorphic forms of any degree.

\begin{theorem} Let $X$ be a projective manifold, $A \subset X$ a submanifold not of general type. Suppose that $N_A$ is big, $\bQ-$effective 
and almost ample (or that $N_A$ is $\bQ-$effective, almost ample and one of the conditions of (3.8) or (3.11) holds). Then 
\begin{enumerate}
\item The image of the canonical map
$$ \pi_1(A) \to \pi_1(X) $$
has finite index in $\pi_1(X)$; the index is at most the number of connected components of the general fiber of $A$ over the
rational quotient $W.$ \\
\item $ h^0(\Omega^q_X) \leq h^0(\Omega^q_A) $ for all $q \geq 1.$
\end{enumerate}

\end{theorem} 

\begin{proof} (1) We are using the followig basic fact: if $f: X \to Y$ is a surjective holomorphic map between normal compact complex spaces,
then the image of $f_*: \pi_1(X) \to \pi_1(Y)$ has finite index in $\pi_1(Y)$ [Ca91,1.3], where the index is also computed. 
Since $X$ is smooth, then the same still holds for meromorphic 
dominant $f$
(blow up $X$ to make $f$ holomorphic; this does not change $\pi_1$). \\
We apply this to the rational quotient $f: X \dasharrow W$ to conclude via 3.6 and 3.7 that the image of 
$$ \pi_1(A) \to \pi_1(W) $$
has finite index. Now the claim follows from the isomorphism
$$\pi_1(X) \simeq \pi_1(W) $$ 
(Koll\'ar [Ko93,5.2]). \\
(2) This follows immediately by considering a holomorphic model of $f \vert A$ and by observing $h^0(\Omega^q_X) = h^0(\Omega^q_W),$ the fibers of $f$
being rationally connected. 
\end{proof}

We are now applying (4.1)  to special cases.

\begin{corollary} Let $A \subset X$ be an abelian variety embedded in a projective manifold $X$. Suppose that $N_A$ is almost ample, $\bQ-$effective and big.
\begin{enumerate}
\item $q(X) \leq \dim A$ and the Albanese map is surjective. 
\item The image of
$$ \pi_1(A) \to \pi_1(X) $$
has finite index in $\pi_1(X)$ and therefore $\pi_1(X)$ is almost abelian. 

\end{enumerate}
\end{corollary}

\begin{proof} (1) is a direct consequence of the almost ampleness of $N_A$ while
(2) is clear from (4.1).
\end{proof} 

\begin{corollary} If $A \subset X$ is a simply connected submanifold not of general type with almost ample, $\bQ-$effective and big normal bundle,
then $\pi_1(X) $ is finite. 
\end{corollary} 

\begin{remark} {\rm In (4.1) - (4.3) the assumption on almost ampleness is only used to make sure that $A$ is not contained in the
indeterminacy locus of the rational quotient. Thus we can omit ``almost ampleness'' by requiring that $A$ is not contained in the
indeterminacy locus. }
\end{remark} 

If $A \subset X$ is a submanifold whose normal bundle is almost ample without assumption on the Kodaira dimension or, more generally, 
generically ample (possible even with the generic spannedness assumption replaced by 
a generic nefness assumption),
then one still would expect that the 
image of $\pi_1(A) \to \pi_1(X)$
has finite index. Here is a partial result in terms of the Shafarevitch map [Ko93],[Ca95]
$$ sh: X \dasharrow {\rm Sh}(X),$$
which is almost holomorphic.  

\begin{proposition} Let $X$ be a projective manifold. Let $A \subset X$ be a submanifold or a local
complete intersection. Suppose that \\
(1) $N_A$ is almost ample,
or that \\
(2) $A$ contains a local complete intersection $B$ such that $N_{B/X}$ is ample. \\
If $\dim A \geq \dim Sh(X),$ resp. $\dim B \geq \dim Sh(X)$, then the image of $\pi_1(A) \to \pi_1(X)$ has finite index.
\end{proposition}

\begin{proof} We use the following fact which was communicated to me by F.Campana: \\
(*) If $sh \vert A$ dominates ${\rm Sh }(X),$ then the image of $\pi_1(A) $ in $\pi_1(X)$ has finite index in $\pi_1(X).$ \\
To prove (*), 
let $h: \tilde X \to X$ be the universal cover and $\tilde A = h^{-1}(A).$ It suffices to show that $\tilde A$ has only
finitely many irreducible components. Consider a general fiber $F$ of $sh$ and let $\tilde F = h^{-1}(F).$ Then $\tilde F$ is compact. 
Now every irreducible component $\tilde A_i$ of $\tilde A$ dominates $A,$ hence $\tilde A_i$ meets $F.$ Since $\tilde A \cap \tilde F$
has only finitely many components, we conclude (*). \\
Now an application of (3.7) yields (1). \\  
(2) By (3.7) $sh \vert B$ dominates ${\rm Sh(X)},$ hence $sh \vert A$ dominates $\rm {Sh(X)},$ and we can apply (*).  

\end{proof}

Of course, this gives nothing when $\dim X = \dim Sh(X),$ i.e. when there is no compact subvariety of positive dimension through the very
general point of the universal cover $\tilde X.$ In that case $X$ is of general type or $\chi(\sO_X) = 0$ [CP04].

\setcounter{equation}{0}
\section{A criterion for rational connectivity}

The following well-known conjecture is due to Mumford.

\begin{conjecture} Let $X$ be a projective manifold. Assume that
$$ H^0(X,(\Omega^1_X)^{\otimes m}) = 0 $$
for all $m \in \bN.$ Then $X$ is rationally connected. 
\end{conjecture} 

This is slightly weaker than the conjecture that $\kappa^+(X) = - \infty$ implies rational connectedness (see (2.2)). 
In fact, $ H^0(X,(\Omega^1_X)^{\otimes m}) = 0 $ implies that
$$ H^0(X,\Gamma(\Omega^1_X)) = 0$$
for all tensor representations $\Gamma$ and in particular $\kappa^+(X) = - \infty.$ 
In this direction we prove

\begin{theorem} Let $X$ be a projective manifold and $L$ be a big line bundle on $X.$
If 
$$ H^0(X,((\Omega^1_X)^{\otimes m} \otimes L)^{\otimes N} ) = 0  \eqno (*)$$
for all $m,N \gg 0$, then $X$ is rationally connected.
\end{theorem}

\begin{proof} The condition (*) implies in particular that
$$H^0(N(mK_X+L)) = 0$$
for $m,N \gg 0.$ Hence $K_X$ is not pseudo-effective, because otherwise $mK_X+L$ would be big. Therefore $X$ is uniruled by [BDPP04]. 
Let $f: X \dasharrow W$ be
the rational quotient ($W$ smooth) and suppose $\dim W > 0.$ 
By eventually blowing up $X$, we clearly may assume that $f$ is holomorphic. Choose a big line bundle $L'$ on $W$ and choose
$k$ so large that $kL - f^*(L')$ is still big. By substituting $L'$ by a multiple and $k$ by a multiple, we may assume that $kL - f^*(L')$
has a section so that 
$$ kL = f^*(L') + A$$
with $A$ effective. 
By (*) we also have 
$$ H^0(X,((\Omega^1_X)^{\otimes m} \otimes kL)^{\otimes N} ) = 0 $$
for large $m,N.$ Therefore 
$$ H^0(W,((\Omega^1_W)^{\otimes m} \otimes L')^{\otimes N}) \subset H^0(((\Omega^1_X)^{\otimes m} \otimes kL)^{\otimes N} ) = 0 $$
for $m,N \gg 0$ Hence by induction $W$ is rationally connected, contradiction. 
Actually it is sufficient to notice that
$$ H^0(N(mK_W+kL')) = 0$$ 
which shows that $K_W$ is not pseudo-effective, hence $W$ is uniruled, contradiction. 
\end{proof} 

Of course, if $X$ is rationally connected, then
$$ H^0(X,((\Omega^1_X)^{\otimes m} \otimes L)^{\otimes N} ) = 0$$
for $m \gg 0$ and $N \in \bN.$ Just restrict to rational curves on which $T_X$ is ample. 

\begin{remark} {\rm By the same method one shows easily the following. If for every $p$ and every $\sF  \subset \Omega^p_X$
one has $$ H^0(N(m\det \sF +L)) = 0$$
for some fixed big line bundle and every $m,N,$ i.e. if $\det \sF$ is not pseudo-effective,
then $X$ is rationally connected. }
\end{remark}

\begin{corollary} Let $X$ be a projective manifold, $C \subset X$ a (possibly singular) curve.
If $T_X \vert C$ is ample, then $X$ is rationally connected. 
\end{corollary} 

\begin{proof} 
We are going to verify (*) in (5.2). Fix a big line bundle $L$ and choose $m_0$ so large that
$$ (\Omega^1_X)^{\otimes {m_0}} \otimes L \vert C $$
is negative. Set
$$ \sE = ((\Omega^1_X)^{\otimes m} \otimes L)^{\otimes N}$$
with $m \geq m_0$ and $N \in \bN.$ Let $\hat C$ be the formal completion of $X$ along $C.$
Then 
$$H^0(\sE) \subset H^0(\sE \vert \hat C) \subset \bigoplus_{k \geq 0} H^0(\sE \otimes \sI^k/\sI^{k+1}).$$
Here $\sI$ denotes the ideal sheaf of $C.$ 
Since $\Omega^1_X \vert C$ is negative, the sheaves $(\sI^k/\sI^{k+1})/{\rm torsion}$ are negative; in fact, the canonical map
$$ \sI^k/\sI^{k+1} \to S^k(\Omega^1_X \vert C) $$
is generically injective. \\ 
Hence any section in
$$ H^0(C,\sI^k/\sI^{k+1} \otimes \sE) $$
is a torsion section, supported at most in the singular locus of $C.$ Thus if we take $s \in H^0(\sE),$ then the 
restriction to the $k-th$ infinitesimal neighborhood is generically $0.$ Since this holds for all $k,$ we conclude that 
$s = 0.$ This shows (*).
\end{proof}

Similarly we prove:

\begin{proposition} 
Let $X$ be a projective manifold  and suppose that $T_X \vert C$ is nef for some (possibly singular) curve $C$. \\
(1) If $-K_X \cdot C > 0,$ then $X$ is uniruled; \\
(2) if $-K_X \cdot C = 0,$ then $\kappa (X) \leq \dim X - 1.$ 
\end{proposition}

\begin{proof} We adapt the proof of (5.4) and substitute $\Omega^1_X$ by $K_X.$ We obtain for all $m$ and all line bundles $L$:
$$ h^0(mK_X-L) \leq h^0(mK_X-L \vert \hat C) \leq \sum_k h^0(mK_X-L \otimes (\sI^k/\sI^{k+1})/{\rm tor}).$$ 
For $L = \sO$ we obtain (1), and by plugging in some ample $L$, we deduce that $h^0(mK_X-L) = 0$ for all $m$ so that 
$K_X$ cannot be big. 
\end{proof} 

\begin{remark} {\rm Assume that 
$$ H^0((\Omega^1_X)^{\otimes m}) = 0  \eqno (*)$$
for all positive $m.$ Suppose that $K_X$ is pseudo-effective. 
{\it If $K_X$ carries a metric with algebraic singularities [DPS01,2.14]}, then
either $\chi(\sO_X) = 0$ - which is excluded by (*) - or 
$$ H^0(X,\Omega^q_X \otimes mK_X) \ne 0 $$
for some $q$ and infinitely many $m.$ This also contradicts (*). 
Consequently if (*) holds, then $K_X$ cannot have a metric with algebraic singularities of non-negative curvature. 
In general, without the assumption on the metric (but still assuming $K_X$ to be pseudo-effective, of course), one can only conclude that 
$$ \chi(X,\sO((m+1)K_X \otimes \sI(h_m)) = 0$$ 
for all $m \geq m_0$ and all (singular) metrics $h_m$ on $mK_X$ with non-negative curvature (current). 
Actually all cohomology groups vanish:
$$ H^0(X,\sO((m+1)K_X) \otimes \sI(h_m)) = 0.$$ 
}
\end{remark}

\vskip .2cm 
\vskip .2cm \vskip .2cm 
\vskip .2cm 
\vskip .5cm

\begin{tabular}{lcl}
Thomas Peternell \\
Mathematisches Institut \\
Universit\" at Bayreuth\\
D-95440 Bayreuth, Germany \\
thomas.peternell@uni-bayreuth.de\\
FAX: + 921-552785
\end{tabular}

\end{document}